\documentclass[10pt, reqno]{amsart}

\usepackage{amsmath} 
\usepackage{amsthm}
\usepackage{amssymb}
\usepackage{amsfonts}
\usepackage{mathdots}
\usepackage{latexsym}
\usepackage{graphicx} 
\usepackage{setspace} 
\usepackage{hyperref}
\usepackage{cite} 
\usepackage[usenames,dvipsnames]{xcolor}
\usepackage{sidecap}
\usepackage{graphicx}
\usepackage{caption}
\usepackage{subcaption}
\usepackage{flafter}
\usepackage{enumitem}
\usepackage{fancyref}

\usepackage{ amssymb }
\usepackage{multicol}
\usepackage{ upgreek }

\usepackage[sc]{mathpazo}
\usepackage[T1]{fontenc}

\numberwithin{equation}{section}
\theoremstyle{plain}

\theoremstyle{definition}

\begin{document}

 \title[Wetzel's Problem, Paul Erd\H{o}s, and the Continuum Hypothesis]{Wetzel's Problem, Paul Erd\H{o}s, and the Continuum Hypothesis: a mathematical mystery}
 \author{Stephan Ramon Garcia}
     \address{   Department of Mathematics\\
            Pomona College\\
            Claremont, California\\
            91711 \\ USA}
    \email{Stephan.Garcia@pomona.edu}
    \urladdr{\url{http://pages.pomona.edu/~sg064747}}
 \author{Amy L. Shoemaker}
 \thanks{Partially supported by NSF Grant DMS-1265973.}

\maketitle 

	We are concerned here with the curious history of \emph{Wetzel's problem}:
	If $\{f_\alpha\}$ is a family of distinct analytic functions (on some fixed domain) such that for each 
	$z$ the set of values $\{f_\alpha (z)\}$ is countable, is the family itself countable?
	
\begin{center}
		\begin{figure}[h]
			\boxed{\includegraphics[scale=.6]{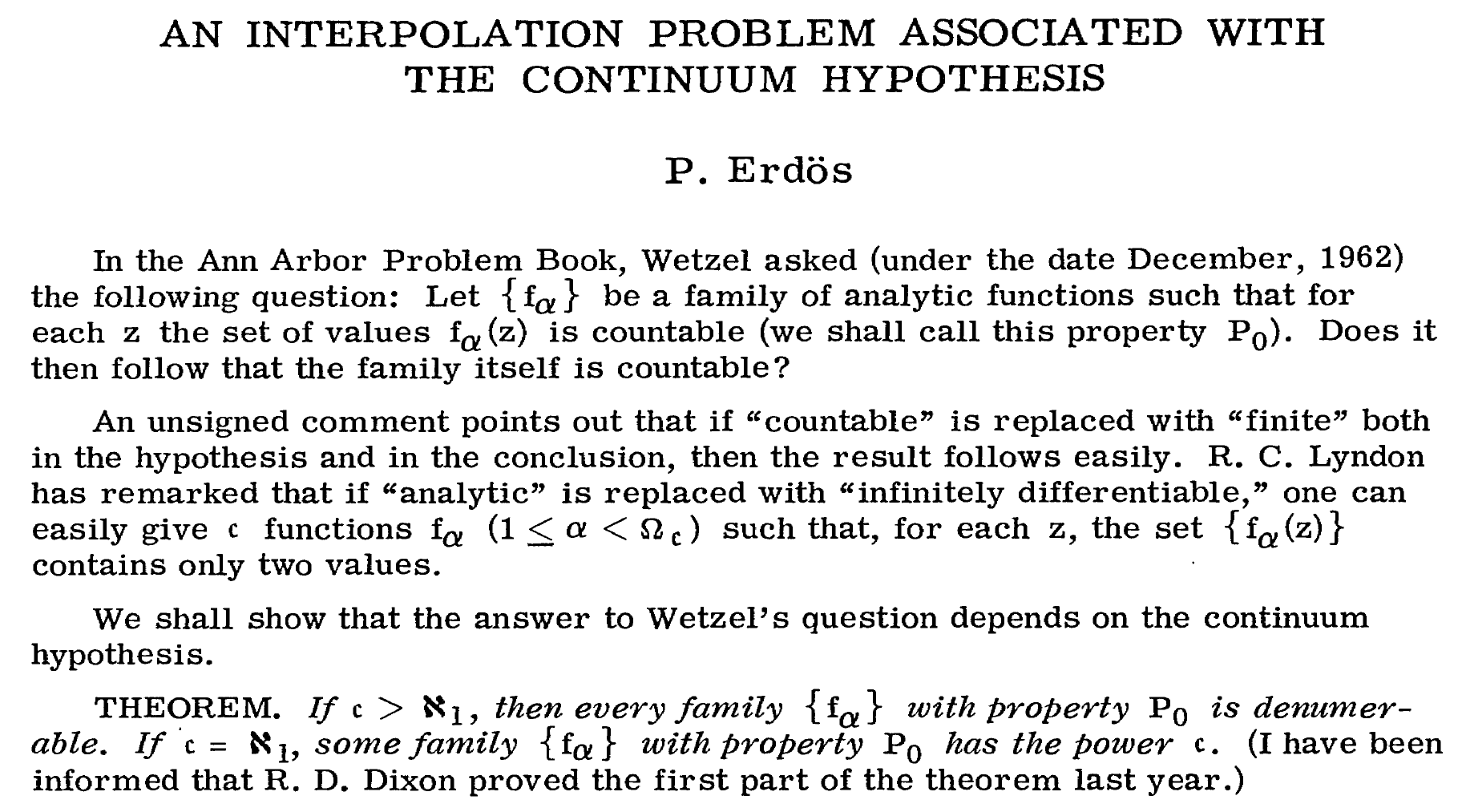}}
			\caption{\footnotesize The first few paragraphs of Erd\H{o}s' paper \cite{erdos}.}\label{erdorig}
		\end{figure}
	\end{center}

	In September 1963, Paul Erd\H{o}s submitted to the Michigan Mathematical Journal a stunning solution to Wetzel's problem (Figure \ref{erdorig}).
	He proved that an affirmative answer is equivalent to the negation of the Continuum Hypothesis.
	Erd\H{o}s ends in an understated manner: ``Paul Cohen's recent proof of the independence of 
	the continuum hypothesis gives this problem some added interest.''
	Together, these results render Wetzel's problem undecidable in ZFC. 

	Erd\H{o}s had a knack for solving ``innocent-looking problems whose solutions shed light on the shape of the mathematical landscape'' 
	\cite[p.~2]{erdosbio1}.   In this case,
	the landscape he revealed was one of underground tunnels, surprising links, and glittering mysteries. 
	However, our interest lies not with the solution itself, but rather with the
	story of how Erd\H{o}s encountered Wetzel's problem in the first place.

	Our first exposure to Wetzel's problem was in \emph{Proofs from the Book} by Aigner and Ziegler.
	``Paul Erd\H{o}s liked to talk about The Book,'' they write, ``in which God maintains the perfect proofs for mathematical theorems, 
	following the dictum of G.H.~Hardy that there is no permanent place for ugly mathematics. Erd\H{o}s also said that you need not believe in God but, 
	as a mathematician, you should believe in The Book'' \cite{thebook}. 
	Erd\H{o}s asked Aigner and Ziegler to assemble a moderate approximation of The Book; 
	included in it was Erd\H{o}s' answer to Wetzel's problem \cite[p.~102-6]{thebook}. 
	
	Regarding the origins of the problem,
	Erd\H{o}s simply asserted that Wetzel posed the question in the Ann Arbor Problem Book in December, 1962.
	Ziegler suggested several mathematicians who might be the Wetzel in question, eventually putting us in
	contact with John E.~Wetzel (Professor Emeritus at the University of Illinois, Urbana-Champaign), who confirmed that the problem was indeed his.
	
\begin{figure}
		\centering
	        \begin{subfigure}[b]{0.25\textwidth}
		        \centering
		   \includegraphics[height=1.3in]{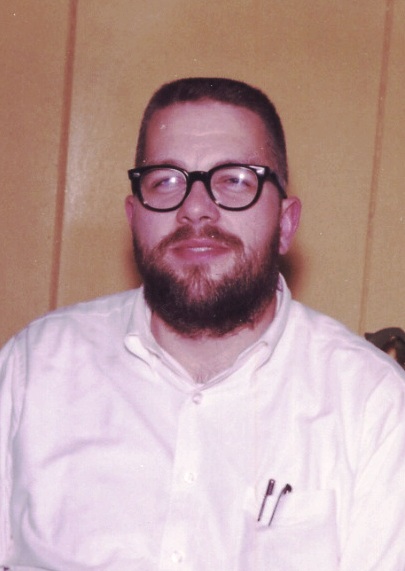}
	                \caption{1963}
	        \end{subfigure} 
	        \begin{subfigure}[b]{0.25\textwidth}\centering
	                \includegraphics[height=1.3in]{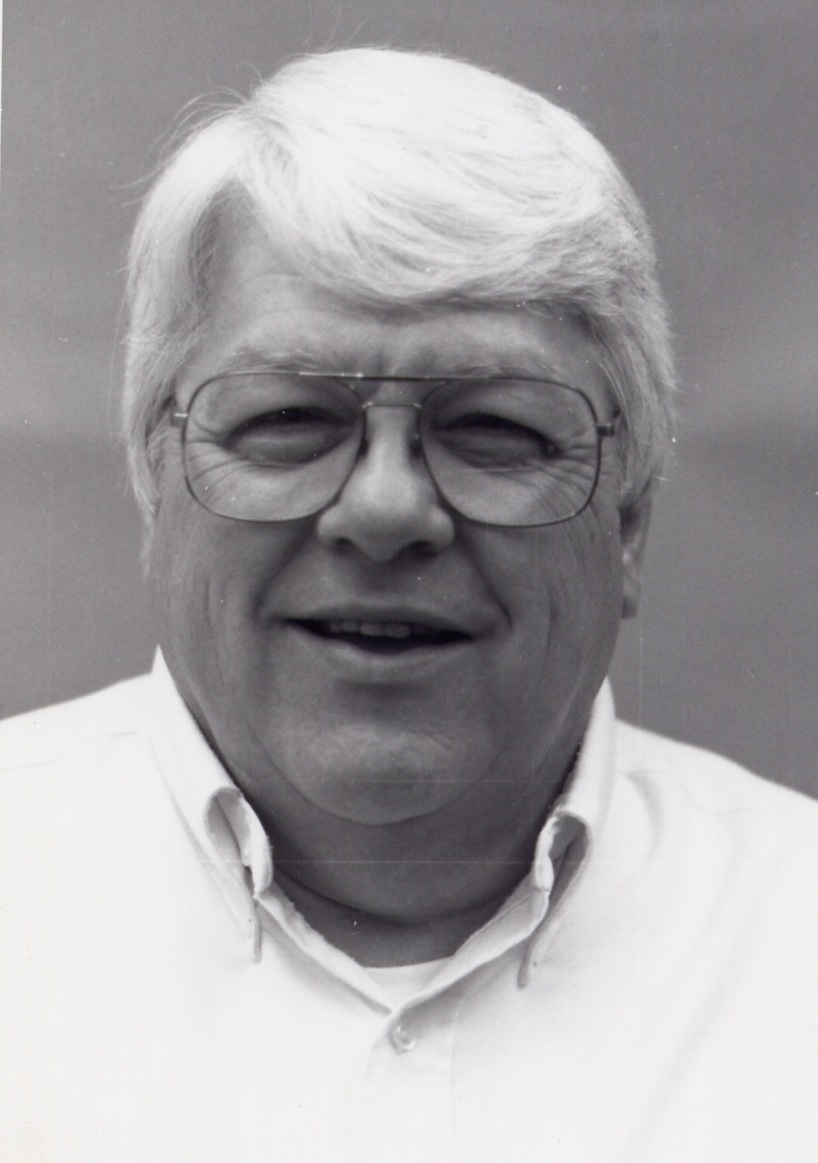}
	                \caption{2008}
	        \end{subfigure}
	        \caption{\footnotesize John E.~Wetzel shortly after he left Stanford (A) and in recent years (B).}
	        \label{wetzel1}
	\end{figure}

	John Wetzel, born on March 6, 1932 in Hammond, Indiana, 
	earned a B.S.~in mathematics and physics from Purdue University in 1954 and went on to study mathematics at Stanford University
	(see Figure \ref{wetzel1}). 
	While studying spaces of harmonic functions on
	Riemann surfaces under Halsey Royden, he posed the following question in his dissertation:
	\begin{quote}\label{wq}\small
		\emph{Let $V$ be a collection of harmonic functions on a Riemann surface $R$ such that for each point $p$ of $R$ the set
		$V_p=\{ v(p): v\in V\}$ is countable. Must $V$ then be countable?} \cite[p.~98]{wetzelthesis}
	\end{quote}
	With regards to the origin of the question, Wetzel explained that Royden had asked him to investigate a specific conjecture:
	\begin{quote}\small
		\emph{I thought about it for a while and eventually showed that what he conjectured was, in fact, not true. I remember reporting to him, thinking that all I had to do was write my work up in good form and I'd be finished with my dissertation; and I remember clearly Royden's reaction that my result would make up perhaps a third of an acceptable dissertation. The question might have had its genesis in the subsequent confusion.} \cite{emailcorr}
	\end{quote}

	Wetzel left Stanford in 1961 to become an Instructor at UIUC, having not yet finished his dissertation.
	He married Rebecca Sprunger in September, 1962 and completed the writing of his dissertation in 1963. 
	During this time, Paul Erd\H{o}s visited the University of Illinois with his mother, who often
	accompanied him as he traveled from campus to campus \cite[p.~7]{erdosbio1}. 

	During one of his casual contacts with Erd\H{o}s, Wetzel recounts,
	\begin{quote}\small
		\emph{I mentioned the question to him, rather timidly, if memory serves.	
		He thought about it briefly and said it was interesting -- and that was the extent of my mathematical contact with him.	
		I don't think he ever told me that he had settled the matter, but every time he visited the campus in the next few years 
		he always asked me if I had any new interesting questions} \cite{emailcorr}.
	\end{quote}
	Upon learning (probably not before 1966, Wetzel said) 
	about Erd\H{o}s' proof, Wetzel wrote to Royden that
	``Erd\H{o}s has showed that the answer to a question I asked in my dissertation is closely tied to the continuum hypothesis!  
	So once again a natural analysis question has grown horns'' \cite{emailcorr}. 
	
	The Ann Arbor Problem Book that Erd\H{o}s mentions seemed most likely to be the Math Club Book, which has its own fiery history. 
	In the words of Peter Duren (Professor Emeritus at the University of Michigan):
	\begin{quote}\small
		\emph{When I came to Ann Arbor in 1962, I learned of a local tradition called ``Math Club,'' an informal gathering 
		of faculty and graduate students which met in the evening every month.  
		A speaker was announced in advance, but the main attraction (in my view) was the series of 3-minute talks, 
		unannounced and often spontaneous, that preceded the announced lecture.  
		There people would tell their colleagues about neat mathematics they had come across, or raise questions, 
		whatever they thought would be of interest.  Afterwards each speaker was invited to write a short summary of his presentation 
		in a book maintained for that purpose.  The Secretary of the Math Club acted as guardian of the book, and both locals and 
		visitors were invited to look through it.  Unfortunately, the book was lost during the Christmas break of 1962-63, on the streets of Chicago.  
		The man then serving as Secretary of the Math Club had carried the book (or books) with him when he drove to 
		Chicago and had left it in his car overnight.  Someone broke into the car and set it on fire, and the Math Club book 
		was lost (among other items, including the car). The Math Club continued to meet, with a new book, 
		but attendance gradually declined and the meetings were discontinued around 1990, as I recall.  
		What Paul Erd\H{o}s called the Ann Arbor Problem Book must have been the Math Club book.  
		But his reference can't be checked, since the original entries for December 1962 no longer exist} \cite{duren}. 
	\end{quote}
	
	However, Wetzel declared
	``I have never visited the University of Michigan; I've never even been to Ann Arbor'' \cite{emailcorr}. 
	This initially led us to conclude that Erd\H{o}s erred in his citation.
	Given his unique manner of doing mathematics and his myriad colleagues, such a slip-up would be understandable. 
	Duren wrote:
	\begin{quote}\small
		\emph{I vividly recall him asking people (including me) at math conferences, ``Where are you located?'', 
		which was the polite way of asking, ``Who are you?  I know I've met you somewhere.''  
		This was only natural, since he traveled so much and met so many mathematicians.  
		It's easy to imagine that he didn't remember correctly where he had seen the problem} \cite{duren}.
	\end{quote}

	After much investigation, we are now able to trace the sequence of events from the original question sparked by 
	Wetzel's dissertation to the problem treated by Erd\H{o}s in his paper. 

	\begin{figure}[h]
		\centering
		\includegraphics[scale=.2]{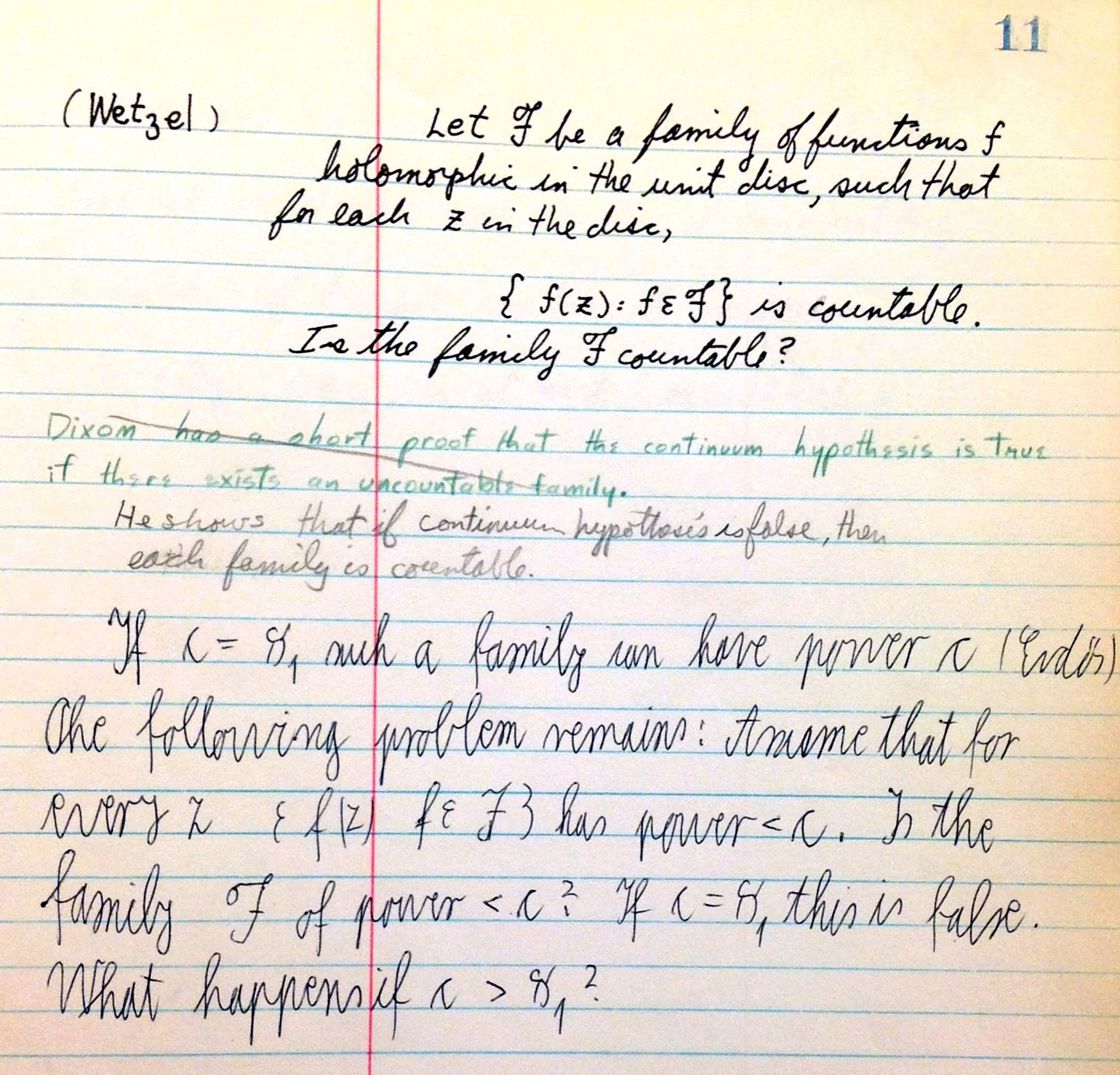}
		\caption{\footnotesize Page eleven of the Boneyard Book (provided by Margaret Lewis \cite{emailcorr}).}\label{boneyard}
	\end{figure} 
	
	Within the UIUC mathematics library is a volume of particular importance for us. Wetzel explained:
	\begin{quote}\small
		\emph{The library maintains a bound volume of blank pages called the Boneyard Book, 
		in which faculty and visiting faculty are encouraged to write problems, 
		including whatever supporting information or commentary they care to include, and faculty looking for interesting problems can browse in it for inspiration }\cite{emailcorr}.
	\end{quote}
	Senior Library Specialist Margaret Lewis recovered the relevant page of the Boneyard Book (see Figure \ref{boneyard}), 
	despite the fact that the volume had been collecting dust in the UIUC archives for decades. 
	Duren recounted that during the early 1960s, 
	Erd\H{o}s traveled frequently between the University of Michigan and the University of Illinois Urbana-Champaign, 
	so conflating the two schools' problem books would have been a natural mistake.

	However, the Boneyard Book led to many more questions. 
	There are clearly four scribes who contributed to the page, so our next task was to identify the mathematicians 
	involved.  While the final entry of the page certainly boasts Erd\H{o}s' distinctive penmanship, 
	the remaining three entries required further scrutiny. 
	
	Wetzel explained that the first entry, though attributed to him, 
	was not written by him: ``I haven't a clue who wrote the problem in our Boneyard Book and wrote my name next to it.  
	I thought for a few moments that it might have been Ranga Rao, with whom I shared an office during my (and his) first year at Illinois, 1961-62.  
	But his handwriting was recognizably `Indian English penmanship,' and I think him unlikely'' \cite{emailcorr}. 
	We therefore decided to investigate the second and third entries before tackling the first.

	Jane Bergman, secretary of the chair of the UIUC math department, 
	suggested that the Dixon mentioned in the Boneyard Book and Erd\H{o}s' article (Figure \ref{erdorig}) was Robert (``Bob'') Dan Dixon:
	\begin{quote}\small
		\emph{I was able to find out that Mr.~Robert Dan Dixon received his Ph.D.~from Ohio State University in 1962 
		in the same year he was hired here as an instructor in Mathematics. In 1964 he was recommended for promotion 
		to assistant professor, but there is no result in his file of that recommendation. I would guess that he moved on 
		at that point, but there is nothing in his file to support my guess. He was born in 1936. I hope you find this helpful.} \cite{emailcorr}
	\end{quote}

	A colleague of Wetzel's, George Robert (``Bob'') Blakley, had more information about Dixon.
	Bob Blakley and Bob Dixon both arrived at UIUC in September 1962 
	as new Ph.D.s and both left in 1964. Blakley, in personal correspondence with the authors and Wetzel, wrote:
	\begin{quote}\small
		\emph{Bob [Dixon] went to the nascent Wright State University in the Dayton area as a founding father.  
		First he founded the math department and headed it.  Then he founded the computer science department and headed 
		it\ldots Late in the last century he retired from WSU, covered with glory\ldots He still manages to fleece me regularly 
		and disreputably in the most varied sorts of bets.  But I think he has given up 100 mile bike rides.}  \cite{blakley}
	\end{quote}
	
	Blakley provided us with three possible ways to contact the elusive Dixon: 
	two email addresses which may or may not have been current, along with his home address. 
	Wetzel sent a letter out to all three addresses, and luckily one route was successful. 
	Dixon responded:
	\begin{quote}\small
		\emph{I was there [UIUC] from the fall of 1962 until the spring of 1964. I remember Erd\H{o}s 
		visiting and may even have had some time with him but I don't remember discussing this problem. 
		The handwriting in the book is puzzling. The entry describing your problem is not familiar to me. 
		The first entry that mentioned me could be by Bob [Blakley] as he generally printed. 
		The second entry refers to me in third person but I could have written it\ldots 
		Although I can't remember any details there is a bit of familiarity. I worked in complex analysis at the 
		time and I had a very interesting course in graduate school that covered the relevant set theory. 
		Lots of problems were thrown around in that group of young faculty.}  \cite{dixon}
	\end{quote}

	Dixon elaborated on ``that group of young faculty'' known as the SixtyTwo Illini Hall Group. 	
	In 1962, UIUC hired twenty new mathematics faculty to add to their faculty of 100; 
	Illini Hall was located across from the mathematics building.
	According to Dixon:
	\begin{quote}\small
		\emph{There was not a sense of privacy about the problems we were investigating. 
		My own work was very specialized and detailed so I had no problems to share. Many
		others, like Jack Wetzel did have problems that they proposed, or pointed us to, that they
		were curious about or needed to solve to get on to the problem they really wanted to 
		work on. We fell then into three overlapping groups, proposers, solvers, and brokers. I 
		was in the solvers group but not particularly successful, Jack may have been in all three. Bob Blakley was in all three but was 
		effective as a broker...There were many others who 
		participated in this interchange but my memory of names is bad. That said it was an 
		experience that had more to do with my career than my own doctoral research. I suspect 
		that was true of several of the other SixtyTwo Illini Hall Group.} \cite{dixon} 
	\end{quote}

	Upon receiving a copy of the Boneyard Book page, Blakley confirmed that he authored the entry in the Boneyard Book that reads, 
	``Dixon has a short proof that the continuum hypothesis is true if there exists an uncountable family.'' 
	Blakley also remarked that he feels ``rather strongly\ldots that Dixon is the third scribe.'' Given Erd\H{o}s' parenthetical remark in his paper that 
	``I have been informed that R.D.~Dixon proved the first part of the theorem last year'' \cite{erdos}, it is likely that Robert Dan Dixon 
	was indeed the scribe of the third entry. Since Dixon expressed that he never spoke to Erd\H{o}s 
	about this problem during their overlapping time at UIUC, the information Erd\H{o}s claimed was relayed to him 
	almost certainly came from the Boneyard Book.
	
	After immersion in the memories and details surrounding the problem, Wetzel recalled:
	\begin{quote}\small
		\emph{I just remembered that I had given a faculty seminar and a departmental colloquium 
		on the substance of my dissertation shortly after arriving at Illinois (even though the dissertation 
		was not yet completely written), and that widened significantly the list of people who might have 
		written the first entry in the Boneyard Book.} \cite{emailcorr}
	\end{quote}
	
	A chance meeting between the first author and John P.~D'Angelo (a Professor at UIUC) finally revealed the most
	likely candidate for the first scribe.  D'Angelo was convinced that the handwriting was that of Lee A.~Rubel:
	\begin{quote}\small
		\emph{Rubel, who died in 1995, often contributed
		to the Boneyard Book. Furthermore, his many interests included
		the interplay between logic and function theory\ldots
		Rubel would have been quite interested in this problem,
		and the handwriting is remarkably similar to that
		of notes he wrote to me around 1979-80.} \cite{dangelo}
	\end{quote}
	Of the possible candidates, Wetzel remarks:
	\begin{quote}\small
		\emph{Lee certainly strikes me as the more likely\ldots 
		I truly don't doubt that Rubel was the author, but I confess that I still find it a little surprising that he never mentioned it to me
		- admitting always the possibility that he did and I have forgotten.} \cite{emailcorr}
	\end{quote}
	The final piece of evidence was a sample of Rubel's handwriting,  
	obtained by D'Angelo from the UIUC archives, which appears to validate this conclusion.	
	In fact, D'Angelo tells us,
	Rubel was the creator of the Boneyard Book.  
	
Let us now return to Erd\H{o}s' paper \cite{erdos}.  It was easy to jump to the conclusion that Erd\H{o}s had erred in his citation, since the problem appears in the Boneyard Book at the University of Illinois, and Erd\H{o}s saw it written there.   
However, there is another scenario that seems to fit the facts more closely.  In the first few lines of his paper (see Figure \ref{erdorig}), Erd\H{o}s cites the ``Ann Arbor Problem Book'' as the source of Wetzel's problem, and mentions ``an unsigned comment'' and a remark by Roger Lyndon (then a professor at Michigan).  Neither of those comments appears on the relevant page of the Boneyard Book (Figure \ref{boneyard}).  

The only explanation is that someone had transported the problem to Michigan and had recorded it either in the Math Club book or in a separate book of open problems.  Peter Duren conjectures that the problem was transferred by Lee Rubel himself, who was making frequent trips to Ann Arbor in those days to work with Allen Shields \cite{duren}.  Erd\H{o}s first saw the problem written there, and came up with the beautiful result presented in [6].  If Erd\H{o}s saw the problem in December 1962, it could well have been in the Math Club book which perished in Chicago.  However, Duren reports having examined the Math Club book for 1962-1991, which now resides in the Bentley Historical Library at UM.  There he found a record that Paul Erd\H{o}s gave a lecture (entitled ``Some Unsolved Problems'') at the Math Club on September 10, 1963, just a week before his paper \cite{erdos} was received by the Michigan Mathematical Journal (on September 18).   Thus it seems far more likely that Erd\H{o}s saw the problem in Ann Arbor during his short visit of 1963, in which case the Ann Arbor Problem Book was not lost in Chicago and may yet be found.  However, Duren has no memory of an independent problem book, nor do his fellow-retirees at Michigan.  In any case, it seems clear that after having submitted his paper \cite{erdos}, Erd\H{o}s saw the problem in the Boneyard Book and learned that Dixon had obtained part of the result independently.

	No story can ever have the entirety of its details pinned down.  As Wetzel said:
	``It may require transfinite induction to bring this matter to a close.'' 
	However, we have identified with a high degree of certainty the trajectory of Wetzel's question as it made its way to Erd\H{o}s. 
	It began in 1961, when Wetzel posed the original question (for harmonic functions on Riemann surfaces) in his evolving dissertation. 	
	When he arrived at UIUC in 1962, he gave a talk on his graduate research.
	Lee Rubel was almost certainly one of the attendees and likely transmitted the problem to Ann Arbor.
	Rubel wrote Wetzel's question in the Boneyard Book in 1962, and Bob Blakley responded with an entry claiming that 
	Bob Dixon had a proof, assuming the truth of the Continuum Hypothesis. Dixon crossed out Blakley's entry and wrote (in third person), 
	``He showed that if the Continuum Hypothesis is false, then each family is countable.'' 
	Dixon's short proof was rediscovered and published by Erd\H{o}s, who went on to prove that an affirmative
	answer to Wetzel's problem is equivalent to the negation of the Continuum Hypothesis.
	Erd\H{o}s' Boneyard Book entry is likely from the fall of 1963, after he had submitted his proof to the
	Michigan Mathematical Journal.

\bibliography{WPECH}{}
\bibliographystyle{plain}

\end{document}